\documentclass[AMA,STIX1COL]{WileyNJD-v2}

\articletype{Article Type}%

\received{19 April 2019}
\revised{-}
\accepted{-}

\begin{document}

\title{The Sequential Test for Chaos}

\author[1]{Marat Akhmet$*$}

\author[2]{Mehmet Onur Fen$\dagger$}

\author[3]{Astrit Tola$\ddagger$}

\authormark{Marat Akhmet \textsc{et al}}

\address[1]{\orgdiv{Department of Mathematics}, \orgname{Middle East Technical University}, \country{Turkey}}

\address[2]{\orgdiv{Department of Mathematics}, \orgname{TED University}, \country{Turkey}}

\address[3]{\orgdiv{Department of Mathematics}, \orgname{Middle East Technical University} \country{Turkey}}

\corres{$*$ \email{marat@metu.edu.tr} \\ $\dagger$ \email{monur.fen@gmail.com} \\ $\ddagger$ \email{astrittola@gmail.com}}

\presentaddress{This is sample for present address text this is sample for present address text}

\abstract[Summary]{This paper reveals a novel numerical method, the sequential test, which approves chaos through sequences of numbers observations. The method alights alongside the Lyapunov exponent and bifurcation diagram test. Explicitly elucidation of the method application for both continuous and discrete systems was given affiliated with the corresponding algorithms. The theoretical results are exemplified on systems satisfying different types of definitions of chaos or numerical methods. The results are supplemented with illustrative graphics.}


\maketitle

\section{Introduction and preliminaries}

At the turn of the nineteenth century, H. Poincar\'{e} started to consider chaotic dynamics. Later on, E. Lorenz \cite{Lorenz63}, Y. Ueda \cite{ueda}, T. Li and J. A. Yorke \cite{Li75} and many others significally developed the theory of chaos.

One can indicate chaos presence either \textit{verifying through a definition} or by \textit{numerical observations}. In the first case, one should consider ingredients of chaos for Devaney's \cite{Dev90}, Li-Yorke's \cite{Li75}, and Poincar\'{e} \cite{Akh17} chaos definitions, and to observe chaos numerically scientists use bifurcation diagrams or Lyapunov exponents \cite{alli, CROB}. Usually numerical observations of chaos does not confirm which type of chaos is satisfied, but in our research, we suggest a novel approach which may approve definition Poincar\'{e} chaos through number sequences observations, that is through \textit{the sequential test}. That is, we give arguments numerically and show that they are for the chaos.

Evaluating the Lyapunov exponent numerical method (LEM) is employed much widely to indicate chaos since it is universally applicable. Nevertheless, the idea is not accepted as a rigorous one since there are examples of nonchaotic systems with positive Lyapunov exponents \cite{alli}. In contrast, the bifurcation diagram analysis method confirms chaos within systems possessing periodic solutions. Our numerical method, the sequential test, is applied to every system which possesses Poisson stable trajectory. Since this method theoretically will request infinitely many iterations to indicate chaos, the same as other methods,  we will apply it and embrace the result for a finite number of iterations.

 Let us denote $\mathbb{N}$ the set of non-negative integers, and consider a metric (X, d) and a map $\pi: \mathbb{T}_+  \times  X \longrightarrow X$, where $\mathbb{T}_+$ is either the set of non-negative real numbers or $\mathbb{N}$, be a semi-flow on X, i.e., $\pi(0, x) = x$ for all $x \in X$, $\pi(t, x)$ is continuous in the pair of variables t and x, and $\pi(t_1,\pi(t_2, x)) = \pi(t_1 + t_2, x)$ for all $t_1, t_2 \in \mathbb{T}_+$, $x \in X$.

A point $x \in X $ is called positively Poisson stable (stable $P^+$) \cite{Sell} if there exists a sequence $\{t_n\}$ satisfying $t_n \longrightarrow \infty$ such that $ \pi (t_n, x) \longrightarrow x$, as $n \longrightarrow \infty$. For a given point $x \in X,$ let $\varTheta_x$ be the closure of the trajectory $T (x) =
\{\pi(t, x) : t \in \mathbb{T}_+\}$. The set $\varTheta_x$ is a quasi-minimal set if the point x is stable $P^+$ and $T(x)$ is contained in a compact subset of X \cite{Sell}.

In paper \cite{Akh21} the definitions of an unpredictable point and Poincar\'{e} chaos was introduced.

\begin{definition} (\cite{Akh21}). 
	A point $p \in X$ and the trajectory through it are unpredictable if there exist a positive number $\varepsilon_0$ (the unpredictability constant) and sequences $\left\{t_n\right\}$ and $\left\{s_n\right\}$, both of which diverge to infinity, such that $f(t_n, p) \to p$ as $n \longrightarrow \infty$ and $d[f(t_n + s_n, p), f(s_n, p)] > \varepsilon_0$ for each $n \in \mathbb{N}$.
\end{definition}

One can see that if a point is unpredictable then it is Poisson stable.

The paper \cite{Akh21} reveals the presence of sensitivity and transitivity in a set $\varTheta_p$ if $p$ is an unpredictable point in X. Their presence in a quasi-minimal set $\varTheta_p$ within $p$ being an unpredictable point, exposed the appearance of chaos, which was a new chaos type named after Poincar\'{e}. Thus the following definition was accepted.

\begin{definition} (\cite{Akh21}). 
	The dynamics on the quasi-minimal set $\varTheta_p$ is called Poincar\'{e} chaotic if p is an unpredictable point.
\end{definition}

In paper \cite{Akh21}, the following Theorem \ref{3} was proven.

\begin{theorem}
	\label{3}
	(\cite{Akh21}). Suppose that $p \in X$ is stable $P^+$ and $T(p)$ is contained in a compact subset of $X$. If $\varTheta_p$ is neither a rest point
	nor a cycle, then it contains an uncountable set of motions everywhere dense and stable $P^+$.
\end{theorem}

It is worth noting that in the paper \cite{Akh21}, it was proved that if $p$ is an unpredictable point, then the dynamics on $p$ is sensitive. That is, there exists a positive number $\widetilde{\varepsilon_0}$ such that for each $x_1 \in \varTheta_p$ and for each positive number $\delta$ there exist a point $x_2 \in \varTheta_p$ and a positive number $\bar{t}$ such that $d(x_1, x_2) < \delta$ and $d(f(\bar{t}, x_1), f(\bar{t}, x_2)) \geq \widetilde{\varepsilon_0}$.

In paper \cite{Akh18}, the definition of unpredictable functions was introduced and it was adapted to the theory of differential equations. In other words, unpredictable functions are considered the solutions of differential equations.
\begin{definition}
	\label{A}
	(\cite{Akh18}).
	A uniformly continuous and bounded function $\vartheta : \mathbb{R} \longrightarrow \mathbb{R}^m$ is unpredictable if there exist positive numbers $\varepsilon_0$, $\delta$ and sequences $\{t_n\}$, $\{s_n\}$ both of which diverge to infinity such that $\left\| \vartheta(t + t_n) - \vartheta(t)\right\| \longrightarrow 0$ as $n\longrightarrow \infty$ uniformly on compact subsets of $\mathbb{R}$ and $\left\|\vartheta(t + t_n) - \vartheta(t)\right\| \geqslant \varepsilon_0$ for each $t \in [s_n - \delta, s_n + \delta]$ and $n \in N$.
\end{definition}

Since to create Poincar\'{e} chaos \cite{Akh21} uniform continuity is not a necessary condition for an unpredictable function $\vartheta$, in paper \cite{Akh18} the Definition \ref{A} was adjusted as follows.

\begin{definition}
	\label{imp3}
	(\cite{Akh18}).
	A continuous and bounded function $\vartheta : \mathbb{R} \longrightarrow \mathbb{R}^m$ is unpredictable if there exist a positive number $\varepsilon_0$ and
	sequences $\{t_n\}$, $\{s_n\}$ both of which diverge to infinity such that $\left\| \vartheta(t + t_n) - \vartheta(t)\right\| \longrightarrow 0$  as $n\longrightarrow \infty$ uniformly on compact subsets of $\mathbb{R}$ and $\left\|\vartheta(t_n + s_n) - \vartheta(s_n)\right\| \geqslant \varepsilon_0$ for each $n \in N$.
\end{definition}

For the convenience of the next discussion, we will call the convergence of the function's shifts on compact subsets and the existance of the sequence $\{t_n\}$ as \textit{Poisson stability of the unpredictable function} or simply Poisson stability, and the existence of the number $\varepsilon_0$ and the sequence $\{s_n\}$ as \textit{unpredictability property} of the function. Thus, a function is unpredictable, if it is \textit{Poisson stable} and admits the \textit{unpredictability property}.

The next, Definition \ref{imp1} and Definition \ref{imp2}, are the instruments for the numerical analysis in this paper. 
\begin{definition}
	\label{imp1}
	(\cite{Akh18}).
	A continuous and bounded function $\vartheta : \mathbb{R} \longrightarrow \mathbb{R}^m$ is unpredictable if there exist a positive number $\varepsilon_0$ and
	sequences $\{t_n\}$, $\{s_n\}$ both of which diverge to infinity such that $\left\| \vartheta(t_n) - \vartheta(0)\right\| \longrightarrow 0$ as $n\longrightarrow \infty$ and $\left\|\vartheta(t_n + s_n) - \vartheta(s_n)\right\| \geqslant \varepsilon_0$ for each $n \in \mathbb{N}$.
\end{definition}

In \cite{Akh18} was also given the definition of the unpredictable sequence. Here, unpredictable sequences are considered as the solutions of discrete equations.

\begin{definition}
	\label{imp2}
	(\cite{Akh18}).
	A bounded sequence ${\kappa_i}$, $i \in \mathbb{N}$, in $\mathbb{R}^p$ is called unpredictable if there exist a positive number $\varepsilon_0$ and sequences $\{\zeta _n\}$, $\{\eta_n\}$, $n \in \mathbb{N}$, of positive integers both of which diverge to infinity such that
	$\left\|\kappa_{\zeta_n} - \kappa_0 \right\| \longrightarrow 0$ as $n\longrightarrow\infty$ and
	$\left\| \kappa_{\zeta_n+\eta_n} - \kappa_{\eta_n} \right\| \geqslant \varepsilon_0$ for each $n \in \mathbb{N}$.
\end{definition}

In course of this definitions, we will suggest the sequential test. Consider the autonomous system of differential equations
\begin{eqnarray} \label{ts_test1}
x'(t)=f(x(t)),
\end{eqnarray}
where $f:\mathbb R^m \to \mathbb R^m$ is a continuous function. Let $x(t)$ be the solution of system (\ref{ts_test1}) with initial condition $x(0)=x_0$, where $x_0$ is a given point in $\mathbb R^m$. 

According to Definition \ref{imp1}, we say that the solution $x(t)$ satisfies the \textit{sequential test}, if it is confirmed numerically that there exist a large natural number $k$ and a positive number $\varepsilon_0$, sequences $\{t_n\}$ and $\{s_n\}$, where $1\leq n\leq k$, for the solution, such that $\left\| x(t_n) - x(0)\right\|=\delta_n$ is a decreasing sequence which becomes close to 0 and the inequality $\left\|x(t_n + s_n)- x(s_n)\right\| > \varepsilon_0$ is valid for every $1 \leq n\leq k$. The largest value of $t_k$ may reach $1.9 \cdot10^6$, as well as the smallest $\delta_k$ are of order $10^{-3}$ For continuous systems. For discrete models the largest values $\zeta_k$ are of order $10^8$ and $\alpha_k$ are small, of order $10^{-6}$. The value of number $k$, as well as the smallness of $\delta_n$, are closely related to the power and facilities offered by a computer, and length of time interval. It is obvious that in different calculations with different $\varepsilon_0$, $\delta_n$ or computers, similar results are not obtained. So the sequences $\{t_n\}$ and $\{s_n\}$ are not unique for a given solution. For convinience, we will call $\{t_n\}$ \textit{the sequence of convergence} and $\{s_n\}$ \textit{the sequence of separation}. In the light of Definition \ref{imp1}, one may say that the solution is unpredictable and the system (\ref{ts_test1}) is Poincar\'{e} chaotic if the sequential test is satisfied. 

It is our hypothesis that if the sequential test works for a sufficient large interval of time, then the Poincar\'{e} chaos is present. In other words, if the test is confirmed for any interval of time preserving the conditions of the test. We suppose that the theorem of this kind can be proved and suggest the assertion as an open problem. The proof of theorem may follow the arguments for Shadowing theorem \cite{ANO,Akh21,PL,ROB,YOR2,PIL}.

In order that, system (\ref{ts_test1}) satisfies the sequential test, we numerically evaluated the sequence $\{t_n\}$ on time interval $(t_{fix},t_{final}]$, where $t_{final}$ is a large real positive number and $0\leq t_{fix}<t_{final}$ is a fixed nunmber. Since $\delta_n$ is a decreasing
sequence which becomes close to 0 then the inequality $\delta_n < \frac{1}{n}$ is valid for some $n$, where $1 \leq n \leq k$. We used this inequality for all systems considered in this paper, on which sequential test is applied. In order to obtain increasing sequences, $\{t_n\}$ and $\{s_n\}$, we set the condition

\textit{(C1) 	\hspace{20pt} $t_{\xi+1}>t_\xi$ and $s_{\xi+1}>s_\xi$, $\xi=1, 2, 3,....$.}\\
Succeeding, we will provide some description of the detailed steps which will be applied later to construct Matlab codes as Algorithm \ref{Alg1}.\\
\\
Let $n=1$ and $\tau_m = m h$, $m=1, 2, ....$.\\
While $\left\| x(\tau_m) - x(0)\right\|<\frac{1}{1}$,\\
fix $t_1=\tau_{m_1}$ the first value of $\tau_m$'s satisfying the inequalities $\left\| x(\tau_m) - x(0)\right\|<1$ and $\tau_{m_1}>t_{fix}$.\\
Set $\tau_i^*=ih$, $i=1, 2, ....$.\\
While $\left\|x(t_1 + \tau_i^*)- x(\tau_i^*)\right\| > \varepsilon_0$,\\
fix $s_1=\tau_{i_1}^*$ the first value of $\tau_i^*$'s satisfying the inequality $\left\|x(t_1 + \tau_i^*)- x(\tau_i^*)\right\| > \varepsilon_0$.\\
\\
Let $n=2$.\\
While $\left\| x(\tau_m) - x(0)\right\|<\frac{1}{2}$,\\
fix $t_2=\tau_{m_2}$ the first value of $\tau_m$'s satisfying the inequalities $\left\| x(\tau_m) - x(0)\right\|<\frac{1}{2}$ and $\tau_{m_2}>t_1$.\\
While $\left\|x(t_2 + \tau_i^*)- x(\tau_i^*)\right\| > \varepsilon_0$,\\
fix $s_2=\tau_{i_2}^*$ the first value of $\tau_i^*$'s satisfying the inequalities $\left\|x(t_2 + \tau_i^*)- x(\tau_i^*)\right\| > \varepsilon_0$ and $\tau_{i_2}^*>s_1$ \\
\\
Let $n=N$, where $1 \leq N \leq k$.\\
While $\left\| x(\tau_m) - x(0)\right\|<\frac{1}{N}$,\\
fix $t_N=\tau_{m_N}$ the first value of $\tau_m$'s satisfying the inequalities $\left\| x(\tau_m) - x(0)\right\|<\frac{1}{N}$ and $\tau_{m_N}>t_{N-1}$.\\
While $\left\|x(t_N + \tau_i^*)- x(\tau_i^*)\right\| > \varepsilon_0$,\\
fix $s_N=\tau_{i_N}^*$ the first value of $\tau_i^*$'s satisfying the inequalities $\left\|x(t_N + \tau_i^*)- x(\tau_i^*)\right\| > \varepsilon_0$ and $\tau_{i_N}^*>s_{N-1}$. 

\newpage

In what follows, the Matlab code based on the above description will be as follows.
\begin{algorithm}
		\caption{Sequential test for system (\ref{ts_test1})}
	\label{Alg1}
\begin{algorithmic}[1]
	
		\State Input $t_{fix}$
		\State Set $l=t_{fix}$ 
		\State Set $q=0$ 
		\State Input $\varepsilon_0$
		\State Input $nspart$ \Comment{number of iterations}
		\State Set $tmin=0$ 
		\State Set $dt=0.01$
		\State Find $tmax=nspart \cdot dt$ 
		\State Input initial condition $x_0$
		\State Find the numerical solution $x(t)$ of system \ref{ts_test1} for the given interval.
		\For{n = 1 : k} 
	    \For{m = 1 : nspart}
		\If{$\left\| x(\tau_m) - x(0)\right\|< \frac{1}{n}$} 
		\If{$l < \tau_m$} 
		\State $l=\tau_m$
		\State A(n)=l \Comment{the matrix A(n) collects $\tau_m$'s, which satisfy lines 13 and 14 for each n}
		\State \textbf{break} \Comment{reckon the first $\tau_m$ replenishing lines 13 and 14 for each n}
		\EndIf
		\EndIf
		\EndFor
		\EndFor
		\For{n = 1 : k} 
		\For{i = 1 : nspart}
		\If{$\left\| x(A(n)+\tau_i^*) - x(\tau_i^*)\right\|> \varepsilon_0$} 
		\If{q < $\tau_i^*$} 
		\State q=$\tau_i^*$ 
		\State B(n)=q \Comment{the matrix B(n) collects $\tau_i^*$'s, which satisfy lines 24 and 25 for each n}
		\State Display matrices ${A(n)}$ and ${B(n)}$
		\State \textbf{break} \Comment{reckon the first $\tau_i^*$ replenishing lines 24 and 25 for each n}
		\EndIf
		\EndIf
		\EndFor
		\EndFor
\end{algorithmic}
\end{algorithm}

The sequential test is also applicable on solutions of discrete systems. For this purpose, cosider the autonomous discrete system
\begin{eqnarray} \label{ts_test2}
x(i+1)=f(x(i)),
\end{eqnarray}
where $f:\mathbb R \to \mathbb R^m$ is a continuous function. Let $x(i)=x_i$ be the solution of system (\ref{ts_test2}) with initial condition $x(0)=x_0$, where $x_0$ is a given point in $\mathbb{R}^m$.

According to Definition \ref{imp2}, we say that the solution $x_i$ satisfies the \textit{sequential test}, if it is confirmed numerically that there exist a large natural number $k$ and a positive number $\varepsilon_0$, increasing sequences of natural numbes $\{\zeta_n\}$ and $\{\eta_n\}$, where $1\leq n\leq k$, for the solution, such that $\left\| x_{\zeta_n} - x_0\right\|=\alpha_n$ is a decreasing sequence which approaches to 0 and the inequality $\left\|x_{\zeta_n + \eta_n}- x_{\eta_n}\right\| > \varepsilon_0$ is valid for every $1 \leq n\leq k$.  Similarly, as in the case of autonomous systems of differential equations, the sequences $\{\zeta_n\}$ and $\{\eta_n\}$ are not unique for a given solution. For convinience, we will call $\{\zeta_n\}$ \textit{the sequence of convergence} and $\{\eta_n\}$ \textit{the sequence of separation}. Setting side by side with Definition \ref{imp2}, one may say that the solution is unpredictable and the system (\ref{ts_test2}) is Poincar\'{e} chaotic if the sequential test is satisfied.

Carrying on, we will present how to create a MATLAB function for the sequential test for autonomous discrete systems. For this reason, consider system (\ref{ts_test2}) with its solution $x_i$ and initial condition $x_0$, where the solution $x_i$ satisfies the sequential test for a positive value $\varepsilon_0$. Based on the following rationalizing, we will erect the upcoming Algorithm \ref{Alg2}.

In order that, system (\ref{ts_test2}) replenishes the sequential test we numerically evaluated the sequence $\{\zeta_n\}$ on interval $(i_{fix},i_{final}]$, where $i_{final}$ is a large natural number and $0\leq t_{fix}<t_{final}$ is a fixed nunmber. Since $\alpha_n$ is a decreasing sequence which becomes close to 0 then the inequality $\alpha_n < \frac{1}{n}$ is valid for some $n$, where $1 \leq n \leq k$. We adopted this inequality for all discrete systems analyzed in this paper, on which sequential test is implemented. In order to obtain increasing sequences, $\{\zeta_n\}$ and $\{\eta_n\}$, we set the condition \\
\textit{(C2)	\hspace{20pt} $\zeta_{\xi+1}>\zeta_k$ and $\eta_{\xi+1}>\eta_k$, $\xi=1, 2, 3,....$.}\\ 
Following, we will implement some explanation of the detailed steps which will be used later to construct Matlab codes as Algorithm \ref{Alg2}.\\
\\
Let $n=1$ and $m=1, 2, ....$.\\
While $\left\| x_m - x_0\right\|<\frac{1}{1}$,\\
fix $\zeta_1=m_1$ the first value of $m$'s satisfying the inequalities $\left\| x_m - x_0\right\|<1$ and $m_1>i_{fix}$.\\
Set $m^*=1, 2, ....$.\\
While $\left\|x_{\zeta_1 + m^*}- x_{m^*}\right\| > \varepsilon_0$,\\
fix $\eta_1=m_1^*$ the first value of $m^*$'s satisfying the inequality $\left\|x_{\zeta_1 + m^*}- x_{m^*}\right\| > \varepsilon_0$.\\
\\
Let $n=2$.\\
While $\left\| x_m - x_0\right\|<\frac{1}{2}$,\\
fix $\zeta_2=m_2$ the first value of $m$'s satisfying the inequalities $\left\| x_m - x_0\right\|<1$ and $m_2>m_1$.\\
While $\left\|x_{\zeta_2 + m^*}- x_{m^*}\right\| > \varepsilon_0$,\\
fix $\eta_2=m_2^*$ the first value of $m^*$'s satisfying the inequality $\left\|x_{\zeta_2 + m^*}- x_{m^*}\right\| > \varepsilon_0$ and $m_2^*>m_1^*$.\\
\\
Let $n=N$, where $1 \leq N \leq k$.\\
While $\left\| x_m - x_0\right\|<\frac{1}{N}$,\\
fix $\zeta_N=m_N$ the first value of $m$'s satisfying the inequalities $\left\| x_m - x_0\right\|<1$ and $m_N>m_{N-1}$.\\
While $\left\|x_{\zeta_N + m^*}- x_{m^*}\right\| > \varepsilon_0$,\\
fix $\eta_N=m_N^*$ the first value of $m^*$'s satisfying the inequality $\left\|x_{\zeta_N + m^*}- x_{m^*}\right\| > \varepsilon_0$ and $m_N^*>m_{N-1}^*$.

\newpage
In what follows, the Matlab code based on the earlier explanation will be as follows.
\begin{algorithm}[ht]
	\caption{Sequential test for system (\ref{ts_test2})}
	\label{Alg2}
	\begin{algorithmic}[1]

		\State Input $i_{fix}$
		\State Set $l=i_{fix}$ 
		\State Set $q=0$ 
		\State Input $\varepsilon_0$
		\State Input $nspart$ (number of iterations)
		\State Input initial condition $x_0$.
		\State Find the numerical solution $x_i$ of system \ref{ts_test2} for the given interval
		\For{n = 1 : k} 
		\For{m = 1 : nspart}
		\If{$\left\| x_m - x_0\right\|< \frac{1}{n}$} 
		\If{$l < m$}
		\State $l=m$
		\State A(n)=l \Comment{the matrix A(n) collects $m$'s, which satisfy lines 10 and 11 for every n}
		\State \textbf{break} \Comment{reckon the first $m$ replenishing lines 10 and 11 for each n}
		\EndIf
		\EndIf
		\EndFor
		\EndFor
		\For{n = 1 : k} 
		\For{$m^*$ = 1 : nspart}
		\If{$\left\| x_{A(n)+m^*} - x_{m^*}\right\|> \varepsilon_0$} 
		\If{q < $m^*$} 
		\State q=$m^*$ 
		\State B(n)=q \Comment{the matrix B(n) collects $m^*$'s, which satisfy lines 21 and 22 for every n}
		\State Display matrices ${A(n)}$ and ${B(n)}$
		\State \textbf{break} \Comment{reckon the first $m^*$ replenishing lines 21 and 22 for each n}
		\EndIf
		\EndIf
		\EndFor
		\EndFor
	\end{algorithmic}
\end{algorithm}

In this paper, we will apply the Algorithm \ref{Alg1} or Algorithm \ref{Alg2} to construct the sequences of convergence, $\{t_n\}$ or $\{\zeta_n\}$, and the sequence of divergence $\{s_n\}$ or $\{\eta_n\}$, respectively, for concrete models that satisfy the sequential test. In other words, the algorithms are the basis of the test. 

\section{Devaney's Chaos Subdued to the Sequential Test}
\begin{figure}[ht]
	\centering
	\includegraphics[height=6.0cm]{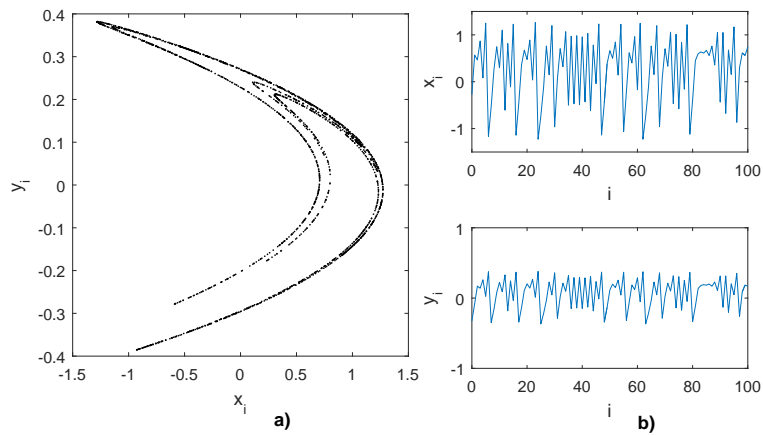}
	\caption{Simulations for the solution of system (\ref{henon}) with the given initial values: \textbf{(a)} the trajectory of the solution, \textbf{(b)} the solution graphs of each coordinate with respect to index i.}
	\label{fig2.3}
\end{figure}

One of the definitions of chaos was provided by Devaney \cite{Dev90} in 1976. To present this definition let us consider the autonomous discrete system 
\begin{eqnarray}
\label{DEV}
x_{i+1}= G(x_i),
\end{eqnarray}
where $G:J \longrightarrow J$, $J$ is the solution space, be continuous. A point $p \in J$ is a \textbf{periodic point} if $G^n(p)=p$, for some $n \geq 1$, and $G^k(p)\neq p$, for $1\leq k <n$. $G:J \longrightarrow J$ is said to be \textbf{topologically transitive} if there is a point $x_0 \in J$ such that the orbit of $x_0$ is dense in $J$. $G:J \longrightarrow J$ is said to have \textbf{sensitive dependence} on initial conditions if there exists $\delta > 0$ such that, for any $x\in J$ and every $\varepsilon>0$, there exists $y \in J$ and $n \geq 0$ such that $|x-y|<\varepsilon$, $|G^n(x)-G^n(y)|> \delta$.
\begin{definition}(\cite{Dev90}).
	The function G is said to be chaotic if:\\
	i. G has sensitive dependence on initial conditions.\\
	ii. G is topologically transitive.\\
	iii. periodic points are dense in J.\\
\end{definition}
One of the most known Devaney chaotic equation is H\'{e}non map. It was introduced on $1975$ by the french astronemer M. H\'{e}non \cite{alli}. Also you can see in \cite{alli} that this map has positive Lyapunov exponents. In his book \cite{Dev90} Devaney proved that H\'{e}non map was Devaney chaotic. Now we will analyse if it is Poincar\'{e} chaotic by using the sequential test. Let us consider the following map:
\begin{eqnarray} \label{henon}
\begin{array}{l}
x_{n+1}=1-1,4x_n^2+y_n\\
y_{n+1}=0.3x_n.
\end{array}
\end{eqnarray}
For this system we took the initial values  $[-0.27518575309954679,-0.32515652033839654]$. Figure \ref{fig2.3}(a) shows the trajectory of system (\ref{henon}) within the initial values and Figure \ref{fig2.3}(b) represents the solution graphs of each coordinate with respect to index i.

We will implement the Sequential Test through Algorithm \ref{Alg2} to system (\ref{henon}) with the fixed initial conditions and $\varepsilon_0 = 2.1$. The index $i$ starts at $0$ and prolongs till $10^8$. As a result, we obtained $4219$ terms for each sequence and selected 11 of them are shown in Table \ref{tab3}.

\begin{table}[ht]
	\centering
	\begin{tabular}{c  c  c  c  c}
		\hline
		n& k& $1/k$ &$\zeta_k$& $\eta_k$\\
		\hline
		$1$& $1$& $1$& $2$& $180$\\
		$2$& $2$& $0.5$& $7$& $181$\\
		$3$& $619$& $0.001616$& $1527239$& $15364$\\
		$4$& $1084$& $0.000923$& $6504212$& $26361$\\
		$5$& $1904$& $0.000525$& $28566825$& $48105$\\
		$6$& $2337$& $0.000428$& $46003272$& $58728$\\
		$7$& $2708$& $0.000369$& $65836673$& $68727$\\
		$8$& $3008$& $0.000332$& $12056560$& $75921$\\
		$9$& $3436$& $0.000291$& $113770045$& $86963$\\
		$10$& $3710$& $0.00027$& $139217264$& $93357$\\
		$11$& $4219$& $0.000237$& $199876783$& $105711$\\
		\hline
	\end{tabular}
	\caption{Selected elements from the sequence of convergence and the sequence of separation obtained by Algorithm \ref{Alg2} applied on system (\ref{henon}).}
	\label{tab3}
\end{table}
\newpage

Next, the results achieved by Algorithm \ref{Alg2} can be displayed graphically. For each element, say $\zeta_\gamma$, within the sequence of convergence $\{\zeta_n\}$, can be drawn a particular graph of solutions of system (\ref{ts_test2}) with initial conditions $x_0$ and $x(\zeta_\gamma)$. We will denote $x_{shift}(i)=x(i+\zeta_\gamma)$ the solution of the system within $x_0=x(\zeta_\gamma)$. In these graphs will be visible the closeness at 0, and the separation bigger than $\varepsilon_0$ between the two solution curves at the corresponding element of the sequence of separation $\{\eta_n\}$, $\eta_\gamma$. We will use this representation on any result obtained by employing Algorithm \ref{Alg2}.

\begin{figure}[ht]
	\centering
	\includegraphics[height=6.0cm]{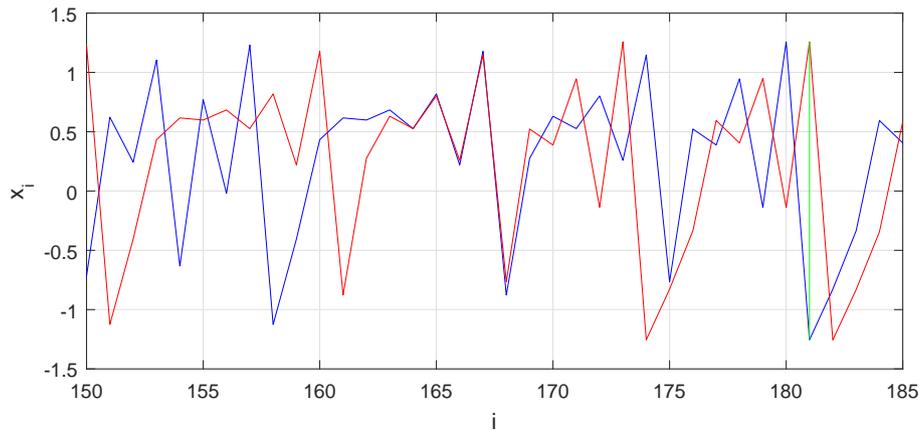}
	\caption{The blue curve shows the graph of solution of system (\ref{henon}), $x(i)$, while the red curve is $x_{shift}(i)$. The green line segment connects the points $(181,x(181))$ and $(181,x_{shift}(181))$ and presents the distance between the solutions at index $i=181$.}
	\label{fig2.4}
\end{figure}

Following, we will represent an individual graph associated with $\zeta_{\gamma}=\zeta_2=7$. Since it is difficult to analize the two dimensional graph, we will show the graph of system (\ref{henon}) for one dimension, the one where the distance between $\omega_i=\omega(i)$ and $\omega_{shift}(i)$, $\omega=x, y$, is bigger than the other dimension at point $i=\eta_\gamma=\eta_{2}=181$. The distance at $\eta_{2}=181$ is bigger in $x$ dimension. In Figure \ref{fig2.4}, the blue curve shows the graph of solution of (\ref{henon}), $x_i=x(i)$, where the initial condition is $X_0 = X(0)$, while the red curve is the solution where the initial value is $X_0 = X(7)$, $x_{shift}(i)=x(i+7)$, where $X(i)=(x(i),y(i))$ and $X_{shift}(i)=(x_{shift}(i),y_{shift}(i))$. The green line segment connects the points $(181,x(181))$ and $(181,x_{shift}(181))$.

\newpage

\begin{figure}[ht]
	\centering
	\includegraphics[height=6.0cm]{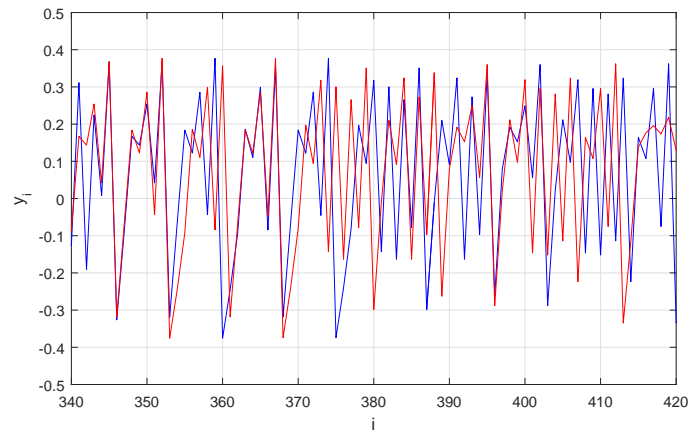}
	\caption{The blue and red curve present the solution of $y(i)$ and $y_{shift}(i)$ of system (\ref{henon}) on the interval [340,420].}
	\label{fig2.4.1}
\end{figure}

The length of the green line segment is $|x(181)-x_{shift}(181)|= 2.5126048968> \varepsilon_0$. In our calculation, we noticed that the distances,\\ $|x(16)-x_{shift}(16)|= 2.4056956028$,\\ $|x(23)-x_{shift}(23)|= 2.2273052694$,\\ $|x(47)-x_{shift}(47)|= 2.3644760759$,\\ $|x(61)-x_{shift}(61)|= 2.2005894834$,\\ $|x(119)-x_{shift}(119)|= 2.267728785$,\\ $|x(126)-x_{shift}(126)|= 2.1601388248$,\\ $|x(174)-x_{shift}(174)|= 2.4024198634$,\\are bigger than $\varepsilon_0$, while the index values are smaller than $\eta_{2}$. Except for these indexes, this is also evident that two-dimensional distances
$\left\|X(40)-X_{shift}(40)\right\|=2.157836143$ and
$\left\|X(72)-X_{shift}(72)\right\|=2.1929991986$. The 2-dimensional distance between the two solution curves at $i=181$ is $\left\|X(181)-X_{shift}(181)\right\|=2.5472810502$. Succeeding, let us consider the closeness of solutions $X(i)$ and $X_{shift}(i)$ of system (\ref{henon}) on the interval [340,420]. Figure \ref{fig2.4.1} presents the graph of solutions $y(i)$ and $y_{shift}(i)$ on [340,420], where the solution curves are nigh.

It is seen from Figure \ref{fig2.4.1}, that the solutions $y(i)$ and $y_{shift}(i)$ are close to each other on the closed intervals [343,353], [361,368], [384,386], [395,400] and [414,417]. The greatest distance between the two solution curves on these intervals is 0.1088336401. If we consider the two-dimensional graph, the solutions $X(i)$ and $X_{shift}(i)$ are close on the closed intervals [344,349], [362,364], [384,385], [393,399] and [414,415]. The greatest two-dimensional distance between the two solution curves on these intervals is 0.1069390621.

\section{Testifying Li-Yorke Chaos}
Li-Yorke chaos was introduced on 1975 in paper \cite{Li75}. To present this definition let $J$ be an interval and consider the autonomous discrete system 
\begin{eqnarray}
x_{i+1}= F(x_i),
\end{eqnarray}
where $F:J \longrightarrow J$ is continuous.
\begin{definition}( \cite{Li75}).
	The function F is said to be chaotic in the sense of Li-Yorke if:\\
	T1. for every $k = 1, 2, 3.....$, there is a periodic point in J having period k.\\
	T2. there is an uncountable set $S \in J$ (containing noperiodic points), which satisfies the following conditions:\\
	A) For every $p, q \in S$ with $p \neq q$\\
	\begin{eqnarray}
	lim_{n \to \infty} sup |F^n(p)-F^n(q)|>0
	\end{eqnarray}
	and 
	\begin{eqnarray}
	lim_{n \to \infty} sup |F^n(p)-F^n(q)|=0.
	\end{eqnarray}
	B) For every $p \in S$ and periodic point $q \in J$, 
	\begin{eqnarray}
	lim_{n \to \infty} sup |F^n(p)-F^n(q)|>0.
	\end{eqnarray}
	
\end{definition}

In their paper \cite{Li75} they also proved that the equation:

\begin{eqnarray} \label{Liyorke}
\begin{array}{l}
x_{i+1}=3.9x_i(1-x_i)
\end{array}
\end{eqnarray}

with initial condition $x(0)=0.5$ is Li-Yorke chaotic.
\begin{figure}[ht]
	\centering
	\includegraphics[height=6.0cm]{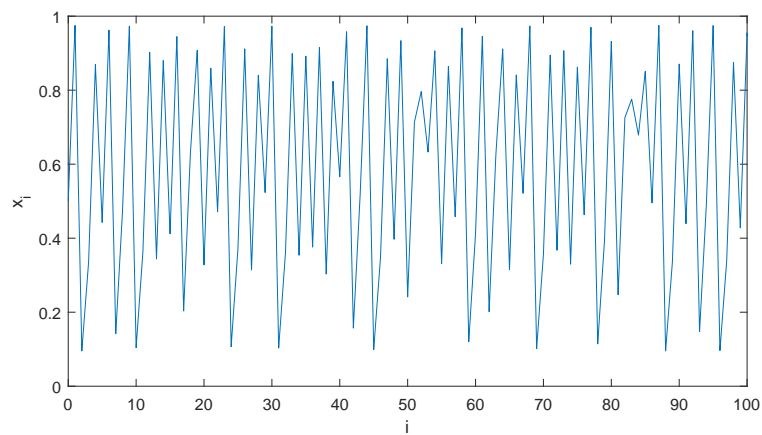}
	\caption{Solution of system (\ref{Liyorke}) with the given initial condition.}
	\label{fig14.1}
\end{figure}

We will apply the sequential test on this system within the set initial condition. The value $i$ starts from $0$ and prolongs till $10^8$. Let $\varepsilon_0 = 0.7$. For this system we obtained $17562$ terms. In Table \ref{tab4} are shown 11 selected elements from the sequence of convergence and the sequence of separation.

\begin{table}[ht]
	\centering
	\begin{tabular}{c  c  c  c  c}
		\hline
		n& k& $1/k$ &$\zeta_k$& $\eta_k$\\
		\hline
		$1$& $1$& $1$& $2$& $2$\\
		$2$& $10$& $0.1$& $40$& $45$\\
		$3$& $1905$& $0.000525$& $1190408$& $15836$\\
		$4$& $3764$& $0.000266$& $4619782$& $30909$\\
		$5$& $5251$& $0.000190$& $9081712$& $42939$\\
		$6$& $6997$& $0.000143$& $16316573$& $57253$\\
		$7$& $8298$& $0.000121$& $22776341$& $67995$\\
		$8$& $9549$& $0.000105$& $29930023$& $78430$\\
		$9$& $11026$& $0.000091$& $39821708$& $90551$\\
		$10$& $12460$& $0.000080$& $50743495$& $10207$\\
		$11$& $17562$& $0.000057$& $99977681$& $16027$\\
		\hline
	\end{tabular}
	\caption{Selected elements from the sequence of convergence and the sequence of separation obtained from Algorithm \ref{Alg2} applied on system (\ref{Liyorke}).}
	\label{tab4}
\end{table}
\newpage
Succeeding, let us graph a particular graph associated with one element within the sequence of convergence presented in Table \ref{tab4}. Let $\zeta_\gamma=\zeta_{10} = 40$. The blue curve shows the graph of solution of (\ref{Liyorke}) where the initial condition is $x_0 = x(0)$, while the red curve is the solution where the initial value is $x_0 = x(40)$. The green line segment connects the points $(45,x(45))$ and $(45,x_{shift}(45))$.

\begin{figure}[ht]
	\centering
	\includegraphics[height=6.0cm]{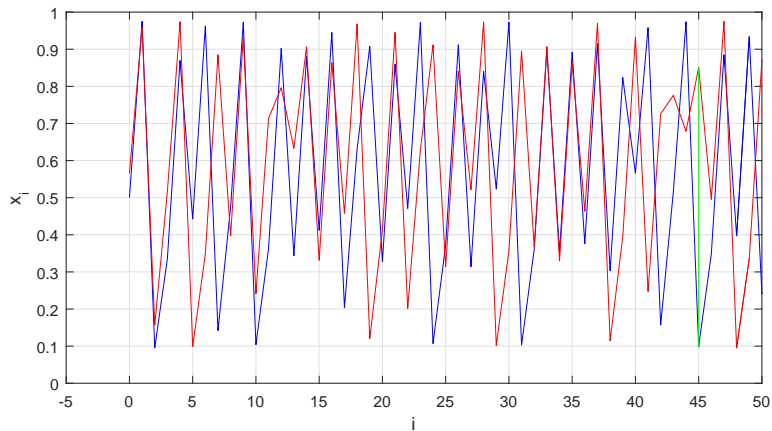}
	\caption{The blue curve shows the graph of solution of the system, $x(i)$, where the initial condition is $x_0 = x(0)$, while the red curve is the solution of system (\ref{Liyorke}) with initial condition $x_0 = x(40)$, $x_{shift}(i)$. The green line segment connects the points $(45,x(45))$ and $(45,x_{shift}(45))$ and presents the distance between the solutions at index $i=45$}
	\label{fig14.2}
\end{figure}

\begin{figure}[ht]
	\centering
	\includegraphics[height=6.0cm]{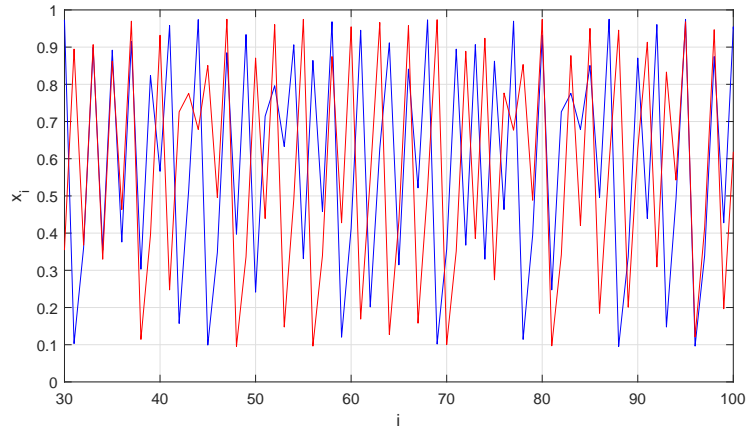}
	\caption{The blue and red curve present the solution of $x(i)$ and $x_{shift}(i)$ of system (\ref{Liyorke}) on the interval [30,100].}
	\label{fig14.2.1}
\end{figure}

The length of the green line segment is $|x(45)-x_{shift}(45)|= 0.7515802112> \varepsilon_0$. In our calculations, we noticed that\\ $|x(7)-x_{shift}(7)|= 0.7429614716$,\\ $|x(19)-x_{shift}(19)|= 0.7870394876$,\\ $|x(24)-x_{shift}(24)|= 0.8045231135$,\\ $|x(31)-x_{shift}(31)|= 0.7915891888$,\\ $|x(41)-x_{shift}(41)|= 0.7105760726$,\\ where the index values are smaller than $\eta_{10}$. Next, let us consider the closeness of solutions $x(i)$ and $x_{shift}(i)$ of the system (\ref{Liyorke}) on the interval [30,100]. Figure \ref{fig14.2.1} displays the graph of solutions on [30,100].

\newpage
It is seen from Figure \ref{fig14.2.1}, that the solutions $x(i)$ and $x_{shift}(i)$ are close to each other on the closed intervals [32,37], [57,58], [65,66], [79,80] and [94,98]. The greatest distance between the two solution curves on these intervals is 0.1187907046.

\section{Bifurcation Diagram Analysis (BDA) and the Sequential Test}

BDA chaotic systems possess periodic solutions \cite{Akh13}. In \cite{Akh16} one can find system (\ref{Period}) and that it is Period-Doubling Route chaotic

\begin{eqnarray} \label{Period}
\begin{array}{l}
x_1'=10(x_2-x_1)\\
x_2'=99.51x_1-x_1x_3-x_2\\
x_3'=x_1x_2-\frac8{3} x_3.
\end{array}
\end{eqnarray}

\begin{figure}[ht]
	\centering
	\includegraphics[height=6.0cm]{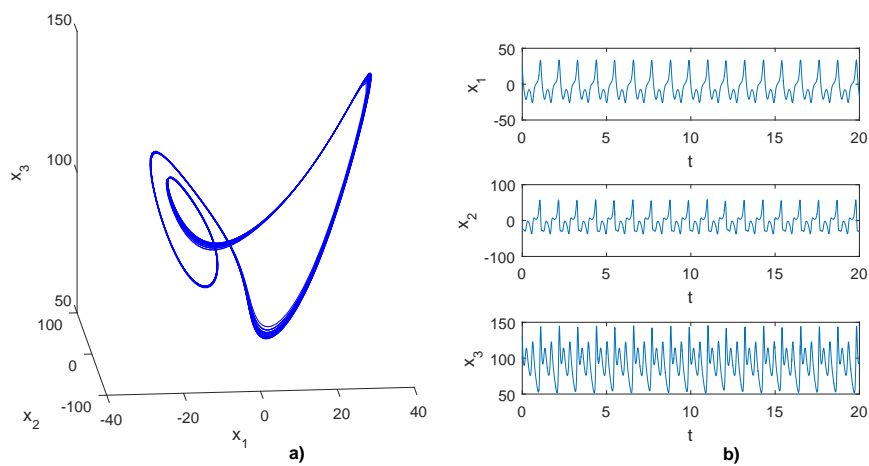}
	\caption{Simulations for the solution of system (\ref{Period}) with the given initial values: \textbf{(a)} the trajectory of the solution, \textbf{(b)} the solution graphs of each coordinate with respect to time t.}
	\label{fig15.1}
\end{figure}

The initial conditions considered are $[23.319088231571342,-15.11725273004282, 130.76383915267931]$.  Figure \ref{fig15.1}(a) presents the trajectory of system (\ref{Period}) within initial conditions while Figure \ref{fig15.1}(b) presents the solution graphs of each coordinate with respect to time t.

We will execute the sequential test through Algorithm \ref{Alg1} to system (\ref{Period}) with the fixed initial conditions. Time interval starts at $0$ and prolongs till $1.3 \cdot 10^6$, partitioned into pieces with distance $0.01$ and $\varepsilon_0 = 45$. In order that the system (\ref{Period}) satisfy the sequential test for $\varepsilon_0=45$, we skipped $t_{fix}=631.36$ while evaluating the sequence of convergence $\{t_n\}$. Within the given conditions and time interval, we found that $k=347$. In Table \ref{tab5}, are shown 10 selected elements from the sequence of convergence and the sequence of separation.

\begin{table}[ht]
	\centering
	\begin{tabular}{c  c  c  c  c}
		\hline
		n& k& $1/k$ &$t_k$& $s_k$\\
		\hline
		$1$& $1$& $1$& $640.35$& $117.87$\\
		$2$& $38$& $0.012195$& $13076.12$& $137.74$\\
		$3$& $79$& $0.006711$& $68393.97$& $160.89$\\
		$4$& $111$& $0.003937$& $118970.29$& $180.69$\\
		$5$& $146$& $0.003311$& $234621.96$& $200.54$\\
		$6$& $175$& $0.002519$& $280456.75$& $217.08$\\
		$7$& $207$& $0.002079$& $331033.07$& $233.63$\\
		$8$& $258$& $0.001706$& $592366.1$& $263.38$\\
		$9$& $310$& $0.001464$& $846071.02$& $293.15$\\
		$10$& $347$& $0.001186$& $1025217.53$& $313$\\
		\hline
	\end{tabular}
	\caption{Selected elements from the sequence of convergence and the sequence of separation obtained from Algorithm \ref{Alg1} applied on system (\ref{Period}).}
	\label{tab5}
\end{table}

Succeeding, the results achieved by Algorithm \ref{Alg1} can be displayed graphically. For each element, say $t_\gamma$, within the sequence of convergence $\{t_n\}$, can be drawn a particular graph of solutions of the system (\ref{ts_test1}) with initial conditions $x_0$ and $x(t_\gamma)$. We will denote $x_{shift}(t)=x(t+t_\gamma)$ the solution of the system within $x_0=x(t_\gamma)$. In these graphs will be visible the closeness at 0, and the separation bigger than $\varepsilon_0$ between the two solution curves at the corresponding element of the sequence of separation $\{s_n\}$, $s_\gamma$. We will use this representation on any result obtained by employing Algorithm \ref{Alg1}.

\begin{figure}[ht]
	\centering
	\includegraphics[height=6.0cm]{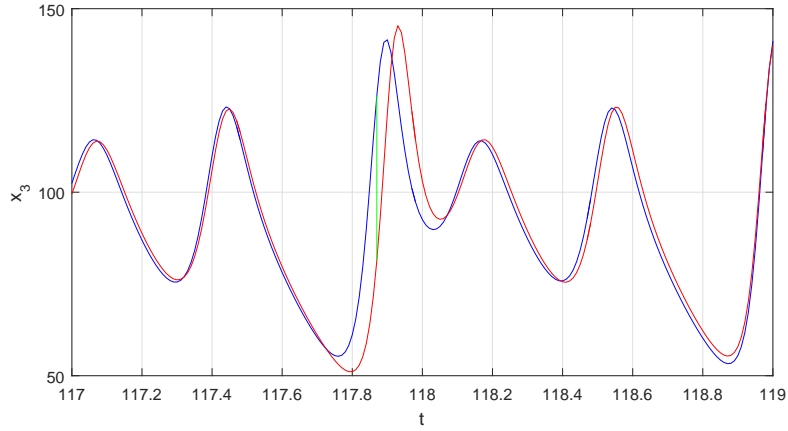}
	\caption{The blue curve shows the graph of solution of the system, $x_3(t)$, while the red curve is $x_{3_{shift}}(t)$. The green line segment connects the points $(117.87,x_3(117.87))$ and $(117.87,x_{3_{shift}}(117.87))$ and presents the distance between the solutions at time $t=117.87$}
	\label{fig15.2}
\end{figure}

Following this description, we will draw a particular graph using $t_\gamma=t_1 = 640.35$. Since it is difficult to analize the three dimensional graph,  we will show the graph of solution for system (\ref{Period}) for one dimension with respect to time, the one where the distance between $x_\omega(t)$ and $x_{\omega_{shift}}(t)$, $\omega=1, 2, 3$, is bigger than the other dimensions at $t=s_\gamma=s_1=117.87$, which is $x_3$ dimension. In Figure \ref{fig15.2}, the blue curve shows the graph of solution of the system (\ref{Period}), $x_3(t)$, where the initial condition is $x_0 = x(0)$, while the red curve is the solution where the initial value is $x_0 = x(640.35)$, $x_{3_{shift}}(t)=x_3(640.35+t)$, where $x(t)=(x_1(t), x_2(t), x_3(t))$ and $x_{shift}(t)=(x_{1_{shift}}(t), x_{2_{shift}}(t), x_{3_{shift}}(t))$. The green line segment connects the points $(117.87,x_3(117.87))$ and $(117.87,x_{3_{shift}}(117.87))$.

\newpage

The length of the green line segment is $|x_3(117.87)-x_{3_{shift}}(117.87)|=44.93479766
$, which is the greatest length between the solution curves until this time in $x_3$ dimension. In $x_1$ dimension the greatest length between the solution curves $|x_1(t)-x_{1_{shift}}(t)|=12.82750911$ occurs at $t=48.51$, while in $x_2$ dimension $|x_2(t)-x_{2_{shift}}(t)|=42.03348827$ occurs at $t=48.47$. The three-dimensional length at $t=117.87$ is $\left\|x(117.87)-x_{shift}(117.87)\right\|=47.03608713$, which is the longest until this time value. 

If someone proceeds calculating further the distances $\left\|x(t)-x_{shift}(t)\right\|$, will notice that all found $s_n$ satisfy the inequality $\left\|x(s_n)-x_{shift}(s_n)\right\|$. Going on, we will show that this is true and for others $t_n$'s. For this purpose we chose four $t_k$ values from Table \ref{tab5}, and in Table \ref{tab21} we presented the lengths $\left\|x(t)-x_{shift_l}(t)\right\|$, where $x_{shift_l}(t)=x(t_l+t)$, for every $s_k$ value of Table \ref{tab5}.

\begin{table}[ht]
	\centering
	\begin{tabular}{c  | c  | c  | c  | c}

			t &$\left\|x(t)-x_{shift_1}(t)\right\|$& $\left\|x(t)-x_{shift_{79}}(t)\right\|$&
			$\left\|x(t)-x_{shift_{175}}(t)\right\|$&
			$\left\|x(t)-x_{shift_{347}}(t)\right\|$\\
			\hline
			$s_1=117.87$	&$47.03608713$ &$47.13469956$ &$47.13803911$ &$47.1349988$  \\
			\hline
			$s_{38}=137.74$	&$48.96987839$ &$49.09002307$ &$49.09409286$ &$49.0903839$  \\
			\hline
			$s_{79}=160.89$ &$48.2635579$ &$48.38466658$ &$48.38876917$ &$48.3850257$  \\
			\hline
			$s_{111}=180.69$&$46.36700321$ &$46.46176646$ &$46.46497564$ &$46.46204434$  \\
			\hline
			$s_{146}=200.54$&$48.85656457$ &$48.96184687$ &$48.96541234$ &$48.96215221$  \\
			\hline
			$s_{175}=217.08$&$49.68721677$ &$49.79847211$ &$49.80224005$ &$49.79879175$  \\
			\hline
			$s_{207}=233.63$&$48.77836421$ &$48.89884399$ &$48.90292512$ &$48.89918692$  \\
			\hline
			$s_{258}=263.38$&$49.70238409$ &$49.81912138$ &$49.82307517$ &$49.81944791$  \\
			\hline
			$s_{310}=293.15$&$48.30769663$ &$48.42876336$ &$48.4328644$ &$48.42909614$  \\
			\hline
			$s_{347}=313$   &$45.86056292$ &$45.98339355$ &$45.98755475$ &$45.98372728$  \\
			
		\end{tabular}	
	\caption{The distances $\left\|x(t)-x_{shift_l}(t)\right\|$ for every $s_k$ value of Table \ref{tab5}, for $l=1, 79, 175, 347$. }
	\label{tab21}
\end{table}

\begin{figure}[ht]
	\centering
	\includegraphics[height=6.0cm]{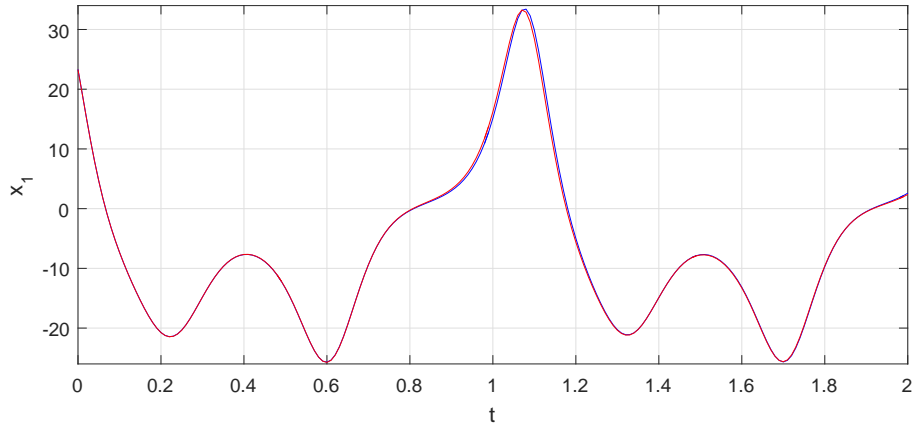}
	\caption{The blue and red curve present the solution of $x_1(t)$ and $x_{1_{shift}}(t)$ of system (\ref{Period}) on the interval [0,2]. Since these solutions curves are close to each other, in the picture they are seen as a single one (red curve).}
	\label{fig15.3}
\end{figure}

Each of the distances shown above is greater than $\varepsilon_0=45$, confirming our claim that \textit{every element of the sequence of separations satisfies the inequality $\left\|x(s_n)-x_{shift_l}(s_n)\right\|>\varepsilon_0$, for every $1\leq l,n \leq k$}. This result is accurate because of periodicity. Possibly, the sequential test can be applied in this way to recognize or to be at least an additional method for analysis of chaos with multiple periods. Following this result, let us consider the closeness of solutions $x(t)$ and $x_{shift}(t)$ of the system (\ref{Period}) on the interval [0,2]. Figure \ref{fig15.3} presents the graph of solutions $x_1(t)$ and $x_{1_{shift}}(t)$ on [0,2], where the solution curves are near.

\newpage

One can notice from Figure \ref{fig15.3}, that the solutions $x_1(t)$ and $x_{1_{shift}}(t)$ are close to each other on the closed intervals [0,0.82], [1.3,1.57], [1.67,1.71] and [1.81,1.94]. The greatest distance between the two solution curves on these intervals is 0.092787609. If we consider the three-dimensional graph, the solutions $x(t)$ and $x_{shift}(t)$ are close on the interval [0.03,0.49]. From our calculations, the greatest three-dimensional distance between the two solution curves on this interval is 0.096167732.

\section{The Sequential Test vs. Lyapunov Exponent Criterium}
Many papers were done by applying LEM, to show that a dynamical system is chaotic. The definition of this numerical method is the following.
\begin{definition}(\cite{alli}).
	Let f be a smooth map on $\mathbb{R}^m$, let $\mathcal{J}_n=Df^n(v_0)$, where $Df^n(v_0)$ denote the first derivative matrix of the $n$th iterate of f, and for $k=1, 2,.....,m$, let $r_n^k$ be the length of the $k$th longest orthogonal axis of the ellipsoid $\mathcal{J}_nU$, where $U$ is a unit sphere, for an orbit with initial point $v_0$ during the first $n$ iterations. The $k$th \textbf{Lyapunov number} of $v_0$ is defined by 
	\begin{eqnarray}
	\begin{array}{l}
	L_k=lim_{n\to \infty}(r_k^n)^{1/n}
	\end{array}
	\end{eqnarray}
	if this limit exists. The $k$th \textbf{Lyapunov exponent} of $v_0$ is $h_k=lnL_k$. The sequence of iterates $x_0, f(x_0), f^2(x_0)=f(f(x_0)),.....,f^n(x_0)...$ defines an \textit{orbit}.
\end{definition}
Following, we apply the sequence test to confirm chaos that has been approved by Lyapunov exponent method.

\subsection{R\"{o}ssler System}
It was firstly introduced by Otto E. R\"{o}ssler in his paper written on $1976$ \cite{rossler76}. In this paper \cite{rossler76}, it was showed that this system is chaoic giving several arguments to show it. The system that we used, which has positive Lyapunov exponents \cite{SPR}, is the following system:
\begin{eqnarray} \label{Rossler}
\begin{array}{l}
x_1'=-x_2-x_3\\
x_2'=x_1+0.2x_2\\
x_3'=0.2+x_1x_3-5.7x_3
\end{array}
\end{eqnarray}

\begin{figure}[ht]
	\centering
	\includegraphics[height=6.0cm]{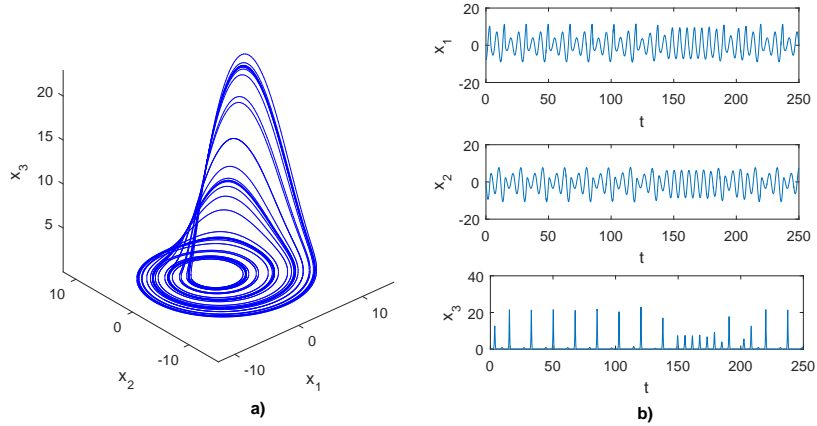}
	\caption{Simulations for the solution of system (\ref{Rossler}) with the given initial values: \textbf{(a)} the trajectory of the solution, \textbf{(b)} the solution graphs of each coordinate with respect to time t.}
	\label{fig 2.1}
\end{figure}

The initial values considered are $[-7.9208550704681606,-0.32213157410506699, 0.01470711076246217]$. Figure \ref{fig 2.1}(a) shows the trajectory of system (\ref{Rossler}) having the fixed initial conditions and Figure\ref{fig 2.1}(b) represents the solution graphs of each coordinate with respect to time t.

\newpage
We will implement the sequential test through Algorithm \ref{Alg1} to system (\ref{Rossler}) with the chosen initial condition on interval \\ $[0, 1.9 \cdot 10^6]$ divided into pieces with distance $0.01$ and $\varepsilon_0 = 22$. In order that the system (\ref{Rossler}) satisfy the sequential test for $\varepsilon_0=22$, we skipped $t_{fix}=35$ while finding the sequence of convergence, $\{t_n\}$. For this system we obtained $157$ terms for each sequence, and 10 of them are shown in Table \ref{tab7}.

\begin{table}[ht]
	\centering
	\begin{tabular}{c  c  c  c  c}
		\hline
		n& k& $1/k$ &$t_k$& $s_k$\\
		\hline
		$1$& $1$& $1$& $64.38$& $20.85$\\
		$2$& $16$& $0.0625$& $2253.47$& $512.43$\\
		$3$& $30$& $0.033333$& $19922.92$& $916.43$\\
		$4$& $56$& $0.017857$& $97107.19$& $2251.43$\\
		$5$& $77$& $0.012987$& $270473.75$& $2796.21$\\
		$6$& $89$& $0.011236$& $405263.82$& $3024.03$\\
		$7$& $111$& $0.009009$& $707522.94$& $3620.64$\\
		$8$& $126$& $0.007937$& $984335.56$& $3754.99$\\
		$9$& $146$& $0.006849$& $1503902.36$& $3848.64$\\
		$10$& $157$& $0.006369$& $1878637.91$& $3848.75$\\
		\hline
	\end{tabular}
	\caption{Selected elements from the sequence of convergence and sequence of separation obtained from Algorithm \ref{Alg1} applied on system (\ref{Rossler}).}
	\label{tab7}
\end{table}

\begin{figure}[ht]
	\centering
	\includegraphics[height=6.0cm]{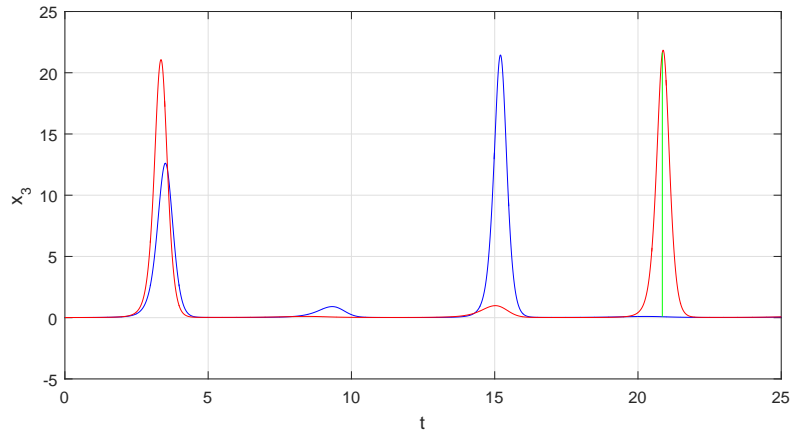}
	\caption{The blue curve shows the graph of solution of the system, $x_3(t)$, while the red curve is $x_{3_{shift}}(t)$. The green line segment connects the points $(20.85,x_2(20.85))$ and $(20.85,x_{2_{shift}}(20.85))$ and presents the distance between the solutions at time $t=20.85$}
	\label{fig 2.2}
\end{figure}

After we acquired the results, let us present the graph of $t_1 = 64.38$, which will be our selected $t_\gamma$. Since it is difficult to analize the three-dimensional graph, we will  show graph of solution for system (\ref{Rossler}) for one dimension with respect to time, the one where the distance between $x_\omega(t)$ and $x_{\omega_{shift}}(t)$, $\omega=1, 2, 3$, is bigger than the other dimensions at point $t=s_\gamma=s_1=20.85$. From our results, the distance at $s_1=20.85$ is bigger in $x_3$ dimension. In Figure \ref{fig 2.2}, the blue curve shows the graph of solution of (\ref{Rossler}), $x_3(t)$, where the initial condition is $x_0 = x(0)$, while the red curve is the solution where the initial value is $x_0 = x(64.38)$, $x_{3_{shift}}(t)=x_3(64.38+t)$, where $x(t)=(x_1(t), x_2(t), x_3(t))$ and $x_{shift}(t)=(x_{1_{shift}}(t), x_{2_{shift}}(t), x_{3_{shift}}(t))$. The green line segment connects the points $(20.85,x_2(20.85))$ and $(20.85,x_{2_{shift}}(20.85))$.

\newpage

If someone would calculate the distance of the green segment would have found that $|x_3(20.85)-x_{3_{shift}}(20.85)|=21,57182863
$, which is the greatest length between the solution curves until this time in $x_3$ dimension. In $x_1$ dimension the greatest length between the solution curves $|x_1(t)-x_{1_{shift}}(t)|=7.887399774$ occurs at $t=20.46$, while in $x_2$ dimension $|x_2(t)-x_{2_{shift}}(t)|=8.506374252$ occurs at $t=19.49$. The three-dimensional length at $t=20.85$ is $\left\|x(20.85)-x_{shift}(20.85)\right\|=22.08237084
$, which is the largest until this time value.

\begin{figure}[ht]
	\centering
	\includegraphics[height=6.0cm]{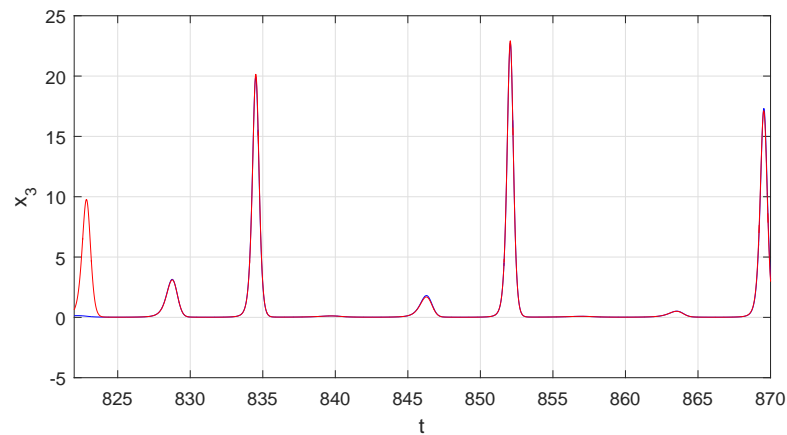}
	\caption{The blue and red curve present the solution of $x_3(t)$ and $x_{3_{shift}}(t)$ of system (\ref{Rossler}) on the interval [822,870].}
	\label{fig 2.3}
\end{figure}

Following this result, let us consider the closeness of solutions $x(t)$ and $x_{shift}(t)$ of the system \ref{Rossler} on the interval [822,870]. Figure \ref{fig 2.3} presents the graph of solutions $x_3(t)$ and $x_{3_{shift}}(t)$ on [822,870], where the solution curves are imminent.

It is seen from Figure \ref{fig 2.3}, that the solutions $x_3(t)$ and $x_{3_{shift}}(t)$ are close to each other on the closed interval [823.79,869.51]. The greatest distance between the two solution curves on this interval is 0.099134999. If we consider the three-dimensional graph, the solutions $x(t)$ and $x_{shift}(t)$ are close on the closed interval [823.8,865.06]. The greatest three-dimensional distance between the two solution curves on this interval is 0.099931352.

\subsection{Ikeda Map}
In this subsection, we considered the following equation taken from the book \cite{alli}, where it was proven that equation (\ref{Ikeda}) has positive Lyapunov exponents.

\begin{eqnarray} \label{Ikeda}
\begin{array}{l}
x_{n+1}=1+0.9(x_ncos(\tau_n)-y_nsin(\tau_n))\\
y_{n+1}=0.9(x_nsin(\tau_n)+y_ncos(\tau_n))\\
\end{array}
\end{eqnarray}
where $\tau_n=0.4-(\frac{6}{1+x_n^2+y_n^2})$

\begin{figure}[ht]
	\centering
	\includegraphics[height=6.0cm]{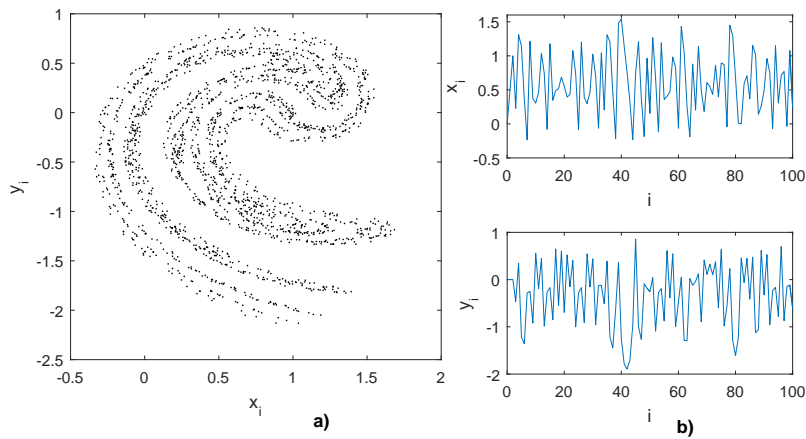}
	\caption{Simulations for the solution of system (\ref{Ikeda}) with the given initial values: \textbf{(a)} the trajectory of the solution, \textbf{(b)} the solution graphs of each coordinate with respect to index i.}
	\label{fig12.1}
\end{figure}

For this system we took as initial conditions $[0,0]$. Figure \ref{fig12.1}(a) shows the trajectory of system (\ref{Ikeda}) within initial conditions and Figure \ref{fig12.1}(b) represents the solution graphs of each coordinate with respect to index i.

We will perform the sequential test through Algorithm \ref{Alg2} to system (\ref{Ikeda}) with the chosen initial condition and $\varepsilon_0 = 2$. Index interval starts from $0$ and prolongs till $3 * 10^6$. As a result, we obtained $754$ terms for each sequence, and  selected 11 of them are shown in Table \ref{tab14}.

\begin{table}[ht]
	\centering
	\begin{tabular}{c  c  c  c  c}
		\hline
		n& k& $1/k$ &$t_k$& $s_k$\\
		\hline
		$1$& $1$& $1$& $2$& $34$\\
		$2$& $16$& $0.0625$& $1193$& $384$\\
		$3$& $95$& $0.010526$& $204042$& $1175$\\
		$4$& $160$& $0.00625$& $1236723$& $2221$\\
		$5$& $201$& $0.004975$& $2394484$& $2920$\\
		$6$& $289$& $0.00346$& $8255283$& $4048$\\
		$7$& $354$& $0.002825$& $14874258$& $4947$\\
		$8$& $407$& $0.002457$& $22484165$& $5858$\\
		$9$& $520$& $0.001923$& $51130431$& $7508$\\
		$10$& $693$& $0.001443$& $161661433$& $9873$\\
		$11$& $754$& $0.001326$& $295347460$& $10735$\\
		\hline
	\end{tabular}
	\caption{Selected elements from the sequence of convergence and the sequence of separation obtained from Algorithm \ref{Alg2} applied on system (\ref{Ikeda}).}
	\label{tab14}
\end{table}

\begin{figure}[ht]
	\centering
	\includegraphics[height=6.0cm]{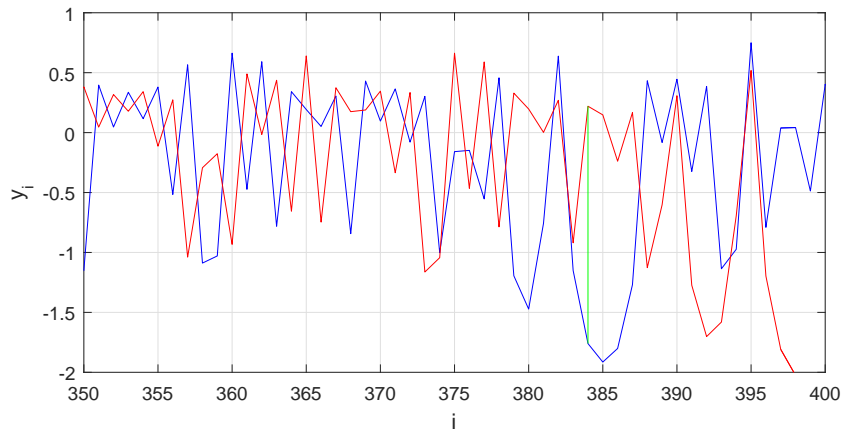}
	\caption{The graph for system (\ref{Ikeda}) with respect to $y_i$ and index $i$. The blue line shows the graph of solution of the system, $y(i)$, while the red line is $y_{shift}(i)$. The green line segment connects the points $(384,y(384))$ and $(384,y_{shift}(384))$ and presents the distance between the solutions at index $i=384$}
	\label{fig12.2}
\end{figure}

After we acquired the results, let us sketch the graph for one of them, say the graph associated with $\zeta_\gamma=\zeta_{16} = 1193$. Since it is difficult to analize the two-dimensional graph, we will show the graph of system (\ref{Ikeda}) for one dimension, the one where the distance between $\omega_i=\omega(i)$ and $\omega_{shift}(i)$, $\omega=x, y$, is bigger than the other dimensions at index $i=\eta_\gamma=\eta_{16}=384$. The distance at $\eta_{16}=384$ is bigger in $y$ dimension. In Figure \ref{fig12.2}, the blue curve shows the graph of solution of (\ref{Ikeda}), $y_i=y(i)$, where the initial condition is $X_0 = X(0)$, while the red curve is the solution where the initial value is $X_0 = X(1193)$, $y_{shift}(i)=y(1193+i)$, where $X(i)=(x(i),y(i))$ and $X_{shift}(i)=(x_{shift}(i),y_{shift}(i))$. The green line segment connects the points $(384,y(384))$ and $(384,y_{shift}(384))$.

\begin{figure}[ht]
	\centering
	\includegraphics[height=6.0cm]{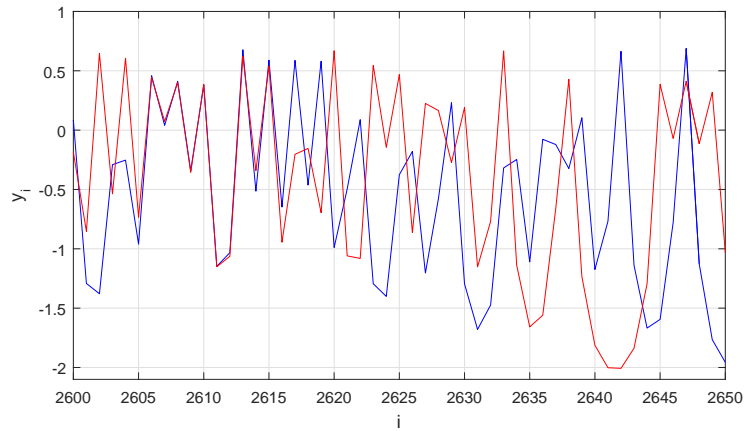}
	\caption{The blue and red curve present the solution of $y(i)$ and $y_{shift}(i)$ of system (\ref{Ikeda}) on the interval [2600,2650].}
	\label{fig12.3}
\end{figure}

From our calculations, we can noticed that the length of the green segment $|y(384)-y_{shift}(384)|=1.980393236<\varepsilon_0$ is not the greatest until the index $i=384$. Someone can find that\\ $|y(22)-y_{shift}(22)|= 2.350274265$,\\ $|y(38)-y_{shift}(38)|= 2.085833749$,\\ $|y(40)-y_{shift}(40)|= 2.261857315$,\\ $|y(42)-y_{shift}(42)|= 2.175361502$,\\ $|y(55)-y_{shift}(55)|= 2.143415044$,\\ $|y(95)-y_{shift}(95)|= 1.983853779$,\\ $|y(106)-y_{shift}(106)|= 2.239737095$,\\ $|y(128)-y_{shift}(128)|= 1.991610248$,\\ $|y(129)-y_{shift}(129)|= 2.243506045$,\\ $|y(131)-y_{shift}(131)|= 2.434072116$,\\ $|y(190)-y_{shift}(190)|= 2.145400124$,\\ $|y(192)-y_{shift}(192)|= 2.212754503$,\\ $|y(212)-y_{shift}(212)|= 2.348997502$,\\ $|y(239)-y_{shift}(239)|= 2.014821546$,\\ $|y(242)-y_{shift}(242)|= 2.167259729$,\\ $|y(282)-y_{shift}(282)|= 2.281796984$,\\ $|y(288)-y_{shift}(288)|= 2.500601744$.\\
Investigating these results, except the distance at $i=95$, all other distances are bigger than $\varepsilon_0$. Except for these indexes, this is also evident at two-dimensional distances\\
$\left\|X(20)-X_{shift}(20)\right\|=2.112599494$,\\
$\left\|X(95)-X_{shift}(95)\right\|=2.007853134$,\\
$\left\|X(120)-X_{shift}(120)\right\|=2.012945907$,\\ $\left\|X(123)-X_{shift}(123)\right\|=2.211804315$,\\ $\left\|X(125)-X_{shift}(125)\right\|=2.147501466$,\\ $\left\|X(126)-X_{shift}(126)\right\|=2.219332021$,\\
$\left\|X(184)-X_{shift}(184)\right\|=2.033164228$,\\
$\left\|X(189)-X_{shift}(189)\right\|=2.084321694$,\\
$\left\|X(239)-X_{shift}(239)\right\|=2.015012267$,\\
$\left\|X(244)-X_{shift}(244)\right\|=2.304006889$,\\
$\left\|X(284)-X_{shift}(284)\right\|=2.151226213$,\\ $\left\|X(285)-X_{shift}(285)\right\|=2.025867605$,\\
The distance for both dimensions at $i=384$ is $\left\|X(384)-X_{shift}(384)\right\|=2.293721422$. Next, let us consider the closeness of solutions $X(i)$ and $X_{shift}(i)$ of the system (\ref{Ikeda}) on the interval [2600,2650]. Figure \ref{fig12.3} presents the graph of solutions $y(i)$ and $y_{shift}(i)$ on [2600,2650], where the solution curves are nigh.

From Figure \ref{fig12.3}, one can notice that the solutions $y(i)$ and $y_{shift}(i)$ are close to each other on the closed interval [2606,2615]. The greatest distance between the two solution curves on this interval is 0.054321108. If we consider the two-dimensional graph, the solutions $X(i)$ and $X_{shift}(i)$ are close on the closed interval [2607,2612]. The greatest two-dimensional distance between the two solution curves on this interval is 0.041486722.

\subsection{Intermittency}
Intermittency or intermittent chaos is a periodic motion where at some specific time chaotic motions burst \cite{Koh84}. The most well-known intermittent system \cite{AI}, which has positive Lyapunov exponents, is:

\begin{eqnarray} 
\label{Intermittency}
\begin{array}{l}
x_1'=10(x_2-x_1)\\
x_2'=-x_1x_3+166.29x_1-x_2\\
x_3'=x_1x_2-\frac8{3} x_3
\end{array}
\end{eqnarray}

\begin{figure}[ht]
	\centering
	\includegraphics[height=6.0cm]{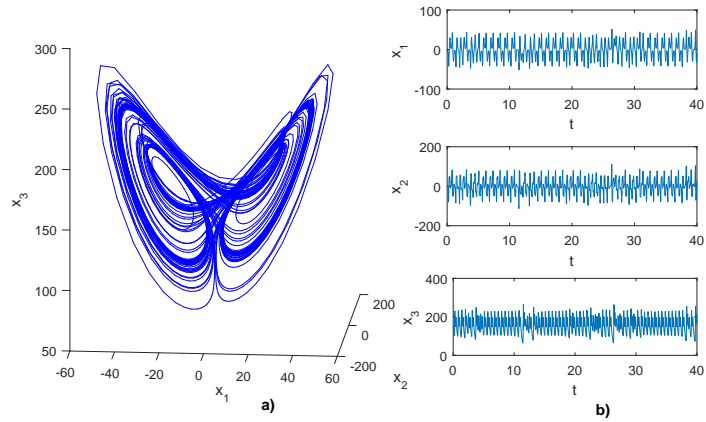}
	\caption{Simulations for the solution of system (\ref{Intermittency}) with the given initial values: \textbf{(a)} the trajectory of the solution, \textbf{(b)} the solution graphs of each coordinate with respect to time t.}
	\label{fig16.1}
\end{figure}

Let the initial conditions be $[-6.9027101537827207, 6.1214285868616205, 146.73481404307805]$. Figure \ref{fig16.1}(a) shows the trajectory of system (\ref{Intermittency}) having the fixed initial conditions, while Figure \ref{fig16.1}(b) represents the solution graphs of each coordinate with respect to time t.

\newpage
We will apply Algorithm \ref{Alg1} on system (\ref{Intermittency}) with the chosen initial condition and $\varepsilon_0 = 200$. Time interval starts from $0$ and prolongs till $1.5 \cdot 10^6$, partitioned into pieces with distance $0.01$. In order that the system (\ref{Intermittency}) satisfy the sequential test for $\varepsilon_0=200$, we skip $t_{fix}=11.64$ while finding the sequence $\{t_n\}$. Within the given conditions and time interval, we found that $k=91$. In Table \ref{tab15}, are shown 10 selected elements from the sequence of convergence and the sequence of separation.

\begin{table}[ht]
	\centering
	\begin{tabular}{c  c  c  c  c}
		\hline
		n& k& $1/k$ &$t_k$& $s_k$\\
		\hline
		$1$& $1$& $1$& $15.43$& $96.98$\\
		$2$& $12$& $0.083333$& $2994.48$& $2322.64$\\
		$3$& $23$& $0.043478$& $12029.93$& $4237.28$\\
		$4$& $31$& $0.032258$& $51983.29$& $6852.58$\\
		$5$& $45$& $0.022222$& $162216.36$& $11141.28$\\
		$6$& $56$& $0.017857$& $401509.63$& $14399.56$\\
		$7$& $63$& $0.015873$& $580535.22$& $16311.43$\\
		$8$& $70$& $0.014286$& $802996.63$& $18649.65$\\
		$9$& $79$& $0.012658$& $1016186.84$& $20581.15$\\
		$10$& $91$& $0.010989$& $1420323.21$& $24758.91$\\
		\hline
	\end{tabular}
	\caption{Selected elements from the sequence of convergence and the sequence of separation obtained from Algorithm \ref{Alg1} applied on system(\ref{Intermittency}).}
	\label{tab15}
\end{table}

\begin{figure}[ht]
	\centering
	\includegraphics[height=6.0cm]{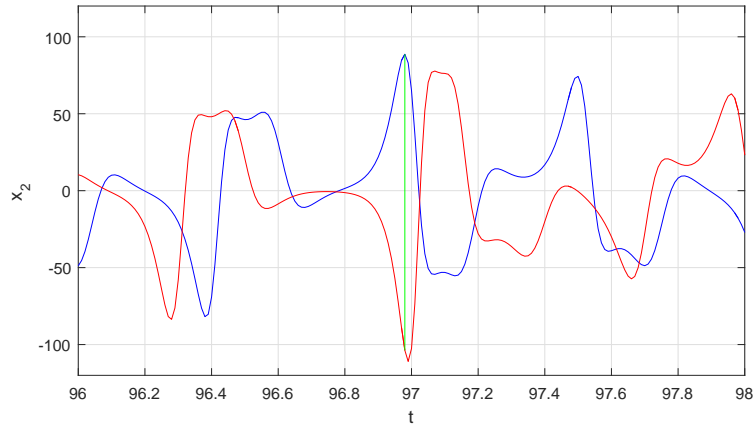}
	\caption{The blue curve shows the graph of solution of system (\ref{Intermittency}), $x_2(t)$, while the red curve is $x_{2_{shift}}(t)$.  The green line segment connects the points $(96.98,x_2(96.98))$ and $(96.98,x_{2_{shift}}(96.98))$ and presents the distance between the solutions at time $t=96.98$}
	\label{fig16.2}
\end{figure}

Succeeding, we will sketch the graph of one element of the sequence of convergence exhibited in Table \ref{tab15}, say for $t_\gamma=t_1 = 15.43$. Since it is difficult to analize the three-dimensional graph, we will show graph of solution for system (\ref{Intermittency}) for one dimension with respect to time, the one where the distance between $x_\omega(t)$ and $x_{\omega_{shift}}(t)$, $\omega=1, 2, 3$, is bigger than the other dimensions at time $t=s_\gamma=s_1=96.98$ , which occurs in $x_2$ dimension. In Figure \ref{fig16.2}, the blue curve shows the graph of solution of (\ref{Intermittency}), $x_2(t)$, where the initial condition is $x_0 = x(0)$, while the red curve is the solution where the initial value is $x_0 = x(15.14)$, $x_{2_{shift}}(t)=x_2(15.43+t)$, where $x(t)=(x_1(t), x_2(t), x_3(t))$ and $x_{shift}(t)=(x_{1_{shift}}(t), x_{2_{shift}}(t), x_{3_{shift}}(t))$. The green line segment connects the points  $(96.98,x_2(96.98))$ and $(96.98,x_{2_{shift}}(96.98))$.

\newpage
The length of the green line segment is $|x_2(96.98)-x_{2_{shift}}(96.98)|=192.44361$, which is the greatest length between the solution curves until this time in $x_2$ dimension. In $x_1$ dimension the greatest length between the solution curves $|x_1(t)-x_{1_{shift}}(t)|=84.04897$ occurs at $t=25.1$, while in $x_3$ dimension $|x_3(t)-x_{3_{shift}}(t)|=173.93722$ occurs at $t=47.51$. The three-dimensional distance at $t=96.98$ is $\left\|x(96.98)-x_{shift}(96.98)\right\|=207.9604947$, which is the biggest until this point. 

\begin{figure}[ht]
	\centering
	\includegraphics[height=6.0cm]{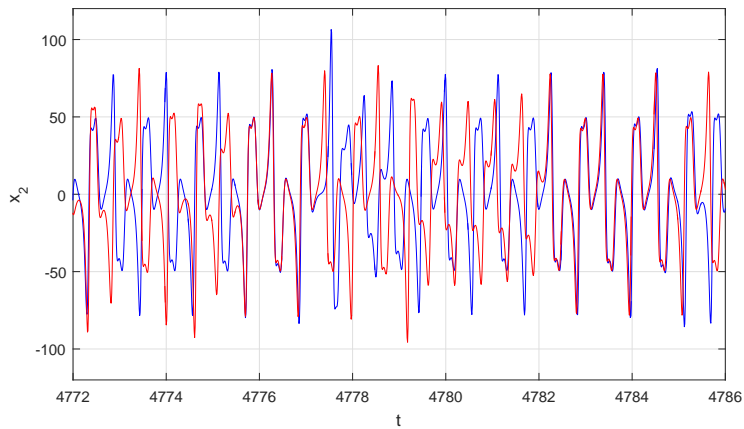}
	\caption{The blue and red curve present the solution of $x_2(t)$ and $x_{2_{shift}}(t)$ of system (\ref{Intermittency}) on the interval [4772,4786].}
	\label{fig16.3}
\end{figure}

Following this result, let us consider the closeness of solutions $x(t)$ and $x_{shift}(t)$ of the system (\ref{Intermittency}) on the interval [4772,4786]. Figure \ref{fig16.3} presents the graph of solutions $x_2(t)$ and $x_{2_{shift}}(t)$ on [4772,4786], where the solution curves are near.

It is seen from Figure \ref{fig16.3}, that the solutions $x_2(t)$ and $x_{2_{shift}}(t)$ are close to each other on the closed intervals [4775.88,4776.07] and [4776.53,4776.56]. The greatest distance between the two solution curves on these intervals is 0.469810895. If we consider the three-dimensional graph, the solutions $x(t)$ and $x_{shift}(t)$ are close on the closed interval [4775.88,4776.07]. The greatest three-dimensional distance between the two solution curves on this interval is 0.527310061.

\end{document}